\documentclass[11pt, a4paper]{article}
\usepackage{geometry}
\usepackage{amsmath, amsthm, amssymb}
\numberwithin{equation}{section}
\begin{document}

\author{Ajai Choudhry}
\title{Circles with four rational points\\ in geometric progression
}
\date{}
\maketitle

\begin{abstract} 
A set of rational points on a curve is said to be in geometric progression if either  the abscissae or the ordinates of the points are in geometric progression. Examples of three points in geometric progression on a circle are already known. In this paper we obtain infinitely many examples of  four points in geometric progression  on a circle with rational radius.
\end{abstract}

\smallskip

\smallskip


\section{Introduction}\label{intro}

Let $C$ be a plane curve defined by an equation $f(x, y)=0$. A set of $n$ rational points on the curve $C$ is said to be in arithmetic or geometric progression if the coordinates $(x_i, y_i), i=1, 2, \ldots, n$, of the $n$ points are such that either the  abscissae $x_i, i=1, 2, \ldots, n$, or the ordinates $y_i, i=1, 2, \ldots, n$, are in arithmetic or in geometric progression, respectively. Several  authors have investigated the existence of sequences of arithmetic or geometric progressions on conics, elliptic curves and hyperelliptic curves (\cite{Br}, \cite{BU}, \cite{Ca}, \cite{Ch1}, \cite{CJ}, \cite{CM1}, \cite{CM2}, \cite{ML}, \cite{Mo1}, \cite{Mo2}, \cite{Ul1}, \cite{Ul2}).

This paper is concerned with geometric progressions on a circle. Çelik, Sadek and  Soydan \cite{CSS} have proved that there exist infinitely many geometric progressions of length 3 on the unit circle $x^2+y^2=1$.  In this paper we obtain infinitely many geometric progressions of length 4 on  circles with rational radii.

\section{Geometric progressions of length 4 on a circle}\label{GPnumex}
We will obtain geometric progressions, with a positive common ratio,  of length 4 on a circle whose centre is on the $x$-axis at $(h, 0)$ and whose radius is $a$, so that the equation of the circle may be written as follows:
\begin{equation}
 (x-h)^2+y^2=a^2, \label{eqcircle}
\end{equation}

If there exist four  points $P_i, i=1, \ldots, 4$, with coordinates $(x_i, y_i)$, $i=1, \ldots, 4$, respectively, on the circle \eqref{eqcircle} such that their abscissae $x_i, i=1, \ldots, 4$, are  in  geometric progression with common positive ratio $r$, then the abscissae of the same points, written in the reverse order, form a geometric progression with common ratio $1/r$. There is, therefore, no loss of generality in assuming that $r > 1$. Since $r$ is positive, all the four abscissae $x_i$ must be  of the same  sign. If these four abscissae are all negative, we note that there is a corresponding circle $(x+h)^2+y^2=a^2$ on which we have four points $(-x_i, y_i), i=1, \ldots, 4$, with positive abscissae and with the same common ratio $r$. We may accordingly assume, without loss of generality, that the abscissae of the  four  points $P_i$ are all positive.

We take both $h$ and $a$ as rational so that  the coordinates of an arbitrary rational point on the circle \eqref{eqcircle} may be written  as $(h+2au/(u^2+1), a(u^2-1)/(u^2+1))$ where $u$ is some rational number. Accordingly, we may take the coordinates of four rational points $P_i, i=1, 2, 3, 4$,  on the circle \eqref{eqcircle} as $(x_i, y_i), i=1, 2, 3, 4$, where
\begin{equation}
x_i=h+2au_i/(u_i^2+1), i=1, \ldots, 4. \label{valx}
\end{equation}

If the abscissae $x_i$ of the four points $P_i, i=1, \ldots, 4$, are in geometric progression with common ratio $r$, we must have $x_{i+1}=rx_i, i=1, 2, 3$, and we thus get the following three conditions, respectively: 
\begin{align}
(u_1^2 + 1)(u_2^2 + 1)(r - 1)h + 2a(ru_1u_2^2 - u_1^2u_2 + ru_1 - u_2)&=0, \label{cond1}\\
(u_2^2 + 1)(u_3^2 + 1)(r - 1)h + 2a(ru_2u_3^2 - u_2^2u_3 + ru_2 - u_3)&=0, \label{cond2} \\
(u_3^2 + 1)(u_4^2 + 1)(r - 1)h + 2a(ru_3u_4^2 - u_3^2u_4 + ru_3 - u_4)&=0. \label{cond3}
\end{align}

The three equations \eqref{cond1}, \eqref{cond2} and \eqref{cond3} may be solved to obtain three values of $h$, namely $h_1, h_2, h_3$ which may be written as follows:
\begin{equation}
\begin{aligned}
h_1&=-2a(ru_1u_2^2 - u_1^2u_2 + ru_1 - u_2)/((u_1^2 + 1)(u_2^2 + 1)(r - 1)), \\
h_2&=-2a(ru_2u_3^2 - u_2^2u_3 + ru_2 - u_3)/((u_2^2 + 1)(u_3^2 + 1)(r - 1)), \\
h_3& =-2a(ru_3u_4^2 - u_3^2u_4 + ru_3 - u_4)/((u_3^2 + 1)(u_4^2 + 1)(r - 1)).
\end{aligned}
\end{equation}

Naturally, if there exists a rational solution of the three simultaneous equations \eqref{cond1}, \eqref{cond2} and \eqref{cond3}, the three values $h_1, h_2, h_3$ must be identical, that is,  we must have $h_1=h_2$ as well as  $h_2=h_3$. The condition $h_1=h_2$ may be written as follows:
\begin{multline}
\quad \quad (ru_2u_3^2 - u_2^2u_3 + u_2u_3^2 + ru_2 + u_2 - u_3)u_1^2 - (u_2^2 + 1)(u_3^2 + 1)ru_1\\
 + ru_2u_3^2 - u_2^2u_3 + u_2u_3^2 + ru_2 + u_2 - u_3=0. \quad \quad \label{condh12}
\end{multline}
Eq. \eqref{condh12} is a quadratic equation in $u_1$,  it does not contain $u_4$,  and it will have a rational solution for $u_1$ if the discriminant $\phi_1(r, u_2, u_3)$ of Eq. \eqref{condh12} with respect to $u_1$ is a perfect square, where 
\begin{multline}
\phi_1(r, u_2, u_3)=(ru_2^2u_3^2 + 2ru_2u_3^2 + ru_2^2 + ru_3^2 - 2u_2^2u_3 + 2u_2u_3^2 + 2ru_2 + r + 2u_2 - 2u_3)\\
\times (ru_2^2u_3^2 - 2ru_2u_3^2 + ru_2^2 + ru_3^2 + 2u_2^2u_3 - 2u_2u_3^2 - 2ru_2 + r - 2u_2 + 2u_3). \label{defphi1}
\end{multline}

Similarly, the condition $h_2=h_3$ yields a quadratic equation in $u_4$ namely,
\begin{multline}
\quad \quad (ru_2^2u_3 - ru_2u_3^2 + u_2^2u_3 - ru_2 + ru_3 + u_3)u_4^2 - (u_2^2 + 1)(u_3^2 + 1)u_4\\
 + ru_2^2u_3 - ru_2u_3^2 + u_2^2u_3 - ru_2 + ru_3 + u_3=0. \quad \quad \label{condh23}
\end{multline}
Eq. \eqref{condh23} does not contain $u_1$ and it  will have a rational solution for $u_4$ if its discriminant $\phi_2(r, u_2, u_3)$  with respect to $u_4$ is a perfect square, where 
\begin{multline}
\phi_2(r, u_2, u_3)=-(2ru_2^2u_3 - 2ru_2u_3^2 + u_2^2u_3^2 + 2u_2^2u_3 - 2ru_2 + 2ru_3 + u_2^2 + u_3^2 + 2u_3 + 1)\\
\times (2ru_2^2u_3 - 2ru_2u_3^2 - u_2^2u_3^2 + 2u_2^2u_3 - 2ru_2 + 2ru_3 - u_2^2 - u_3^2 + 2u_3 - 1). \label{defphi2}
\end{multline}

We performed computer trials to find rational values of $u_2, u_3$ and a positive rational value of $r$ such that both the discriminants 
$\phi_1(r, u_2, u_3)$ and $\phi_2(r, u_2, u_3)$ are perfect squares. We took $u_2=p_1/p_2, u_3=q_1/q_2, r=r_1/r_2$ where $p_i, q_i, r_i, i=1, 2$, are integers such that  $r_1 > r_2 >0$, and performed trials over the range $|p_1| + |p_2| +|q_1| +|q_2| + r_1 +r_2 \leq 100$. This yielded 14 sets of numerical values of $u_2, u_3, r$ such that both $\phi_1(r, u_2, u_3)$ and $\phi_2(r, u_2, u_3)$ are perfect squares, and we thus obtained 14 rational solutions of the simultaneous Eqs. \eqref{condh12} and \eqref{condh23} leading to  14 numerical examples of four points $P_i, i=1, \ldots, 4$, in geometric progression on the  circle \eqref{eqcircle} with a rational radius. In Table \ref{Tnumex} we have listed 6 of these  numerical examples in which  the abscissae of the points $P_i, i=1, \ldots, 4$,  are all less than $10^4$.  
\begin{table}[tbh]
\caption{Geometric progressions of length 4 on a circle}
\label{Tnumex}
\begin{center}
\begin{tabular}{|c||c|c|c|c|c|c|} \hline
S. No. &  $u_2$ &$u_3$ &$r$& $h$ &$a$ & Points $P_1, P_2, P_3, P_4$\\
\hline
1. & 1/3 & 0 & 5/3 & 75  & 50 & (27,14), (45,40), (75,50), (125,0) \\
2. & 3/5 &  1/3 & 8/3 &  447 &  425 & (27,65), (72,200), (192,340), (512,420) \\
3. & 11/23  & 1/5  & 5/2 &  5075 &  4875 & (512, 1716), (1280, 3060), \\
&&&&&& (3200, 4500), (8000, 3900)\\
4. & 11/23 &  4/7  & 5/3  &  -2125 &   3250 & (243, 2226), (405, 2040), \\
&&&&&& (675, 1650), (1125, 0) \\
5. & 11/29 &   1/31  & 16/9 &  4825 & 3367 & (1458, 0), (2592, 2520),\\
&&&&&&  (4608, 3360), (8192, 0) \\
6. & 3/5 &  21/67  & 3 & 2559 &  2465 & (128, 408), (384, 1160),\\
&&&&&& (1152, 2024), (3456, 2296) \\
\hline 
\end{tabular}
\end{center}
\end{table}

\section{An infinitude of geometric progressions}\label{GPinf} We will now use the numerical results of  Section \ref{GPnumex} to  obtain infinitely many circles with  rational radii and having four rational points in geometric progression.
\subsection{Circles in which one or both  points on the $x$-axis are included in the set of four points in geometric progression}\label{xaxiscircle}
We observe from Table \ref{Tnumex} that, in  three of the geometric progressions listed at S. No. 1, 4 and 5, the point $P_4$ lies on the $x$-axis, and the coordinates of $P_4$ may be written as $(h+a, 0)$. On further investigation, we found that there are infinitely many such geometric progressions. While we can find them by solving the simultaneous equations \eqref{condh12} and \eqref{condh23}, we obtain them below by starting {\it ab initio} since this is simpler.

If the point $P_1$ on the circle \eqref{eqcircle} is defined as before, and $P_4$ is taken as $ (h+a, 0)$, and the common ratio of the geometric progression is $r$, then $x_1r^3=h+a$, that is,   $r^3(h + 2au_1/(u_1^2 + 1)) =h+a$, hence we get,
\begin{equation}
h=a(u_1^2-2r^3u_1 + 1)/((u_1^2 + 1)(r^3 - 1)).
\end{equation}

Since the points $P_i, i=1, \ldots, 4$,  are in geometric progression, the abscissae of the points $P_2$ and $P_3$ must be $x_1r$ and $x_2r$ respectively, and hence we get,
\begin{equation}
\begin{aligned}
x_2 &=ar(u_1 - 1)^2/((u_1^2 + 1)(r^3 - 1)), \\
x_3 & = ar^2(u_1 - 1)^2/((u_1^2 + 1)(r^3 - 1)).
\end{aligned}
\label{valx23}
\end{equation}

Since both $P_2$ and $P_3$ lie on the circle \eqref{eqcircle}, we must have $a^2-(x_2-h)^2=y_2^2$ and $a^2-(x_3-h)^2=y_3^2$, and on writing
\begin{align}
y_2 & = v(u_1-1)a/((r^2 + r + 1)(u_1^2 + 1)), \\
y_3 & = w(u_1-1)ra/((r^2 + r + 1)(u_1^2 + 1)),
\end{align}
we get the following two conditions, respectively:
\begin{align}
(u_1 + 1)^2r^4 + 2(u_1 + 1)^2r^3 + (3u_1^2 + 2u_1 + 3)r^2 + (2u_1^2 + 2)r &=v^2, \label{condx2}\\
(u_1 + 1)^2r^2 + (2u_1^2 + 2)r + 2u_1^2 + 2 &= w^2, \label{condx3}
\end{align}
where $v$ and $w$ are some rational numbers.

In the next subsection we obtain  solutions of the simultaneous diophantine equations \eqref{condx2} and \eqref{condx3} that lead  to  geometric progressions containing  
both  points of the circle \eqref{eqcircle} on the $x$-axis, while in the following subsection, we obtain solutions that yield geometric progressions in which only one point of the circle on the $x$-axis is included.

\subsubsection{Geometric progressions containing both points $(h-a, 0)$ and $(h+a, 0)$}\label{bothptsxaxis}

When the geometric progression  begins with the point $P_1 =(h-a, 0)$, we have $u_1=-1$, and now the two conditions \eqref{condx2} and \eqref{condx3} reduce to $4r(r+1)=v^2$ and $4(r+1)=w^2$, respectively. Thus, both $r$ and $r+1$ must be perfect squares, hence we take $r=(t^2-1)^2/(4t^2)$, where $t$ is an arbitrary rational parameter. This leads to the circle  \eqref{eqcircle} where
\begin{equation}
\begin{aligned}
h & = (t^8 - 8t^6 + 30t^4 - 8t^2 + 1)(t^2 + 1)^2, \\
a & = t^{12} - 6t^{10} + 15t^8 - 84t^6 + 15t^4 - 6t^2 + 1,
\end{aligned}
\label{circspl1}
\end{equation}
with the following four points in geometric progression:
\begin{equation}
\begin{aligned}
P_1 & = (128t^6, 0), \\
P_2 &= (32(t - 1)^2(t + 1)^2t^4, 8t^2(t^8 - 6t^6 + 6t^2 - 1)),\\
P_3 & = (8(t - 1)^4(t + 1)^4t^2, 4t(t^{10} - 7t^8 + 6t^6 + 6t^4 - 7t^2 + 1),\\
P_4 & =  (2(t - 1)^6(t + 1)^6, 0).
\end{aligned}
\end{equation}

Taking $t=3$ yields, on appropriate scaling, the numerical example listed at S. No. 5 in Table I.

\subsubsection{Geometric progressions with only one point on the $x$-axis}\label{simgleptxaxis}
We will now solve the simultaneous diophantine equations \eqref{condx2} and \eqref{condx3} to obtain geometric progressions that begin with a point on the circle \eqref{eqcircle} that is not on the $x$-axis and end with the point $(h+a, 0)$ which lies on the $x$-axis.

We write
\begin{equation}
u_1=(1+t)/(1-t),\quad  v=2p/(t-1), \quad w=2q/(t-1), \label{valuvw}
\end{equation}
when Eqs. \eqref{condx2} and \eqref{condx3} may be written as
\begin{align}
r(r+1)(r^2+t^2+r+1)&=p^2, \label{condx2a}\\
(r+1) t^2+r^2+r+1 &=q^2, \label{condx3a}
\end{align}
respectively. 

On eliminating $t$ between Eqs. \eqref{condx2a} and \eqref{condx3a} we get an equation that may be written as follows:
\begin{equation}
(r^2+p+r)(-r^2+p-r)=r(q-r)(q+r). \label{condx23}
\end{equation}
Eq. \eqref{condx23} will be satisfied if there exists a rational number $X$ such that 
\begin{align}
r^2 + p + r &=X(q+r), \label{condx23a}\\
X(-r^2 + p - r) &=r(q-r).  \label{condx23b}
\end{align}

Eqs. \eqref{condx23a} and \eqref{condx23b} are two linear equations in $p$ and $q$, and on solving them, we get
\begin{align}
p &= r(X^2r + X^2 - 2Xr + r^2 + r)/(X^2 - r), \label{valp}\\
q & = -r(X^2 - 2Xr - 2X + r)/(X^2 - r), \label{valq}
\end{align}

Now on substituting the value of $p$ given by \eqref{valp} in Eq. \eqref{condx2a}, we get
\begin{multline}
(r+1)(X^2-r)^2t^2+(r+1)X^4+4r^2(r+1)X^3-2r(2r^3+6r^2+3r+1)X^2\\
+4r^3(r+1)X+r^2(r+1)=0, \quad \quad \label{condx3f}
\end{multline}
and, on writing $t=Y/((r+1)(X^2-r))$, Eq. \eqref{condx3f} may be written as follows:
\begin{align}
Y^2 & =-(r+1)^2X^4-4r^2(r+1)^2X^3+2r(r+1)(2r^3+6r^2+3r+1)X^2 \nonumber \\
 & \quad \quad-4r^3(r+1)^2X-r^2(r+1)^2. \label{ecXY}
\end{align}

Eq. \eqref{ecXY} is a quartic equation in $X$ and $Y$ and it may be considered as a quartic model of an elliptic curve over the function field  $\mathbb{Q}(r)$. If for any numerical rational value of $r$, the rank of the elliptic curve \eqref{ecXY} is positive, we can find infinitely many rational  solutions of Eq. \eqref{ecXY}, and hence also of the simultaneous diophantine equations \eqref{condx2} and \eqref{condx3}, and thus  obtain  geometric progressions on the circle \eqref{eqcircle}  which end with the point $(h+a, 0)$. 

We wrote $r=r_1/r_2$ and performed trials, taking $r_1 > r_2 > 0$, over the range $ r_1+r_2 \leq 20$, and obtained elliptic curves of positive rank when $r=5/3$,  $5/4$ and $ 9/7$. Each of these values of $r$ will yield infinitely many solutions of the simultaneous diophantine equations \eqref{condx2} and \eqref{condx3} and hence we get  infinitely many examples of geometric progressions on the circle  \eqref{eqcircle} including the point $(h+a, 0)$. 

As a numerical example, when $r=5/3$, Eq. \eqref{ecXY} may be written as
\begin{equation}
Y^2 = -64/9X^4 - 6400/81X^3 + 68960/243X^2 - 32000/243X - 1600/81, \label{ecXYex1}
\end{equation}
and the birational transformation defined by 
\begin{equation}
\begin{aligned}
 X& = -3(10x+y - 700 )/(2( 4x -y - 980)), \\
Y &= 98(x^3 - 360x^2 + 14700x - 1000y + 1192000)/(9( 4x -y - 980)^2),
\end{aligned}
\label{biratXY2xy}
\end{equation}
and
\begin{equation}
\begin{aligned}
x & = -(208X^2 - 960X - 81Y + 90)/(2X-3)^2,\\
y & = -18(528X^3 - 920X^2 - 36XY - 1980X - 135Y + 1200)/(2X-3)^3,
\end{aligned}
\label{biratxy2XY}
\end{equation}
reduces Eq. \eqref{ecXYex1} to the Weierstrass form of the elliptic curve given by 
\begin{equation}
y^2 = x^3 - 14700x + 286000. \label{weierec}
\end{equation}

A reference to Cremona's well-known database of elliptic curves \cite{Cr} shows that the rank of the curve \eqref{weierec} is 2 with the  two generators of the Mordell-Weil group being $R_1=(-40, 900)$ and $R_2=(14, 288)$. We can, therefore, find infinitely many rational points on the elliptic curve  \eqref{weierec} using the group law, and then find the corresponding point on the quartic curve \eqref{ecXYex1} using the relations \eqref{biratXY2xy} and thus obtain  the four points $P_i$ in geometric progression on the circle \eqref{eqcircle}.

As a numerical example, the point $R_1$ on the curve \eqref{weierec} yields the four points
\[
(4563, 29666), (7605, 27560), (12675, 22750), (21125, 0),
\]
which are in geometric progression with common ratio $5/3$, on the circle
\begin{equation}
(x+13725)^2+y^2=34850^2.
\end{equation}
Similarly, the point $R_2$ on the curve \eqref{weierec} yields the four points
\[
(88356987, 171988866), (147261645, 198917640), \\
(245436075, 201974850), (409060125, 0)
\]
 which are in geometric progression, with the same common ratio $5/3$, on the circle
\begin{equation}
(x-202590875)^2+y^2=206469250^2.
\end{equation}

Finally we note that there  also exist circles \eqref{eqcircle} which have geometric progressions of four points beginning with the point $(h-a,0)$ which lies  on the $x$-axis while the fourth point does not lie on the $x$-axis. Such examples may be found by following a method similar to that used above for finding geometric progressions ending with the point $(h+a, 0)$. Accordingly, we restrict ourselves to giving only a  numerical example. The circle
\begin{equation}
(x-966123)^2+y^2=963235^2
\end{equation}
has the two points $(2888, 0)$ and $(1929358, 0)$ on the $x$-axis, and it has the following four points in geometric progression with common ratio 8, and beginning with the point $(2888, 0)$:
\[
(2888,0),(23104,196308),(184832,563388),(1478656, 815556).
\]

\subsection{Circles in which the points on the $x$-axis are not included 
in the geometric progression}\label{GPgen}
We will now obtain infinitely many examples of geometric progressions of points $P_1, \ldots$, $P_4$, on the  circle \eqref{eqcircle} such that the points of intersection of the circle with the $x$-axis are not included in the geometric progression.
\subsubsection{Geometric progressions in which the points $P_1$ and $P_4$ are symmetrically situated on either side of the diameter parallel to the $y$-axis}\label{GPsym}

The existence of geometric progressions beginning with the point $(h-a, 0)$ and ending with $(h+a, 0)$ suggests the possibility of geometric progressions in which the two points  $P_1$ and $P_4$ are symmetrically situated on either side of the vertical diameter of the circle. Accordingly, we take the abscissae of the points $P_1$ and $P_4$  as $h - 2au/(u^2 + 1)$ and $h + 2au/(u^2 + 1)$. Since the common ratio is $r$, we have $(h - 2au/(u^2 + 1))r^3=h + 2au/(u^2 + 1)$, and hence we get
\begin{equation}
h=2au(r^3 + 1)/((u^2 + 1)(r^3 - 1)). \label{GPsymvalh}
\end{equation}

Since the abscissae of the points $P_i, i=1, \ldots, 4$, are in geometric progression the abscissae of the points $P_2$ and $P_3$ are given by 
\begin{align}
x_2 &=(h - 2au/(u^2 + 1))r, \label{GPsymvalx2}\\
x_3 &= (h - 2au/(u^2 + 1))r^2, \label{GPsymvalx3}
\end{align}
respectively. 

Since $P_2$ and $P_3$ are rational points on the circle \eqref{eqcircle}, the following two conditions must be satisfied:
\begin{align}
a^2-(x_2-h)^2 &=y_2^2, \label{GPsymcondP2}\\
a^2-(x_3-h)^2 &=y_3^2, \label{GPsymcondP3}
\end{align}
where $y_2$ and $y_3$ are some rational numbers.  Using the values of $h, x_2$ and $ x_3$ given by  \eqref{GPsymvalh}, \eqref{GPsymvalx2}  and \eqref{GPsymvalx3}, respectively, and taking
\begin{align}
y_2& =az_2/((r^2 + r + 1)(u^2 + 1)), \label{GPsymvaly2}\\
y_3 &=az_3/((r^2 + r + 1)(u^2 + 1)),   \label{GPsymvaly3}
\end{align}
the conditions \eqref{GPsymcondP2} and \eqref{GPsymcondP3} may be written as follows:
\begin{align}
 (u^2 - 1)^2r^4 + 2(u^2 - 1)^2r^3 + (u^2 + 3)(3u^2 + 1)r^2 &  \nonumber \\
+ 2(u^4 + 6u^2 + 1)r + (u^2 - 1)^2&=z_2^2, \label{GPsymcondP2a}\\
 (u^2 - 1)^2r^4 + 2(u^4 + 6u^2 + 1)r^3 + (u^2 + 3)(3u^2 + 1)r^2 \nonumber \\
 + 2(u^2 - 1)^2r + (u^2 - 1)^2&=z_3^2, \label{GPsymcondP3a}
\end{align}

When $u=3$ the discriminants, with respect to $r$, of each of the two equations \eqref{GPsymcondP2a} and \eqref{GPsymcondP3a} vanish, and these two equations may now be written as follows: 
\begin{align}
16(r^2 + r + 4)(2r + 1)^2 &=z_2^2, \label{GPsymcondP2b}\\
16(4r^2 + r + 1)(r + 2)^2&=z_3^2. \label{GPsymcondP3b}
\end{align}

By choosing $r$ such that  $r^2 + r + 4$ becomes a perfect square,   we readily obtain the following solution of Eq. \eqref{GPsymcondP2b}:
\begin{equation}
r = (p^2 - 4q^2)/(q(2p + q)), \quad  z_2 = 4(2p^2 + 2pq - 7q^2)(p^2 + pq + 4q^2)/(q^2(2p + q)^2), \label{valrz2}
\end{equation} 
where $p$ and $q$ are arbitrary parameters. 

Now, on using the value of $r$ given by \eqref{valrz2}, and writing
\begin{equation}
p=Xq, \quad z_3=4(X^2 + 4X - 2)Y/(2X + 1)^2, \label{valpz3}
\end{equation}
Eq. \eqref{GPsymcondP3b} may be written as follows:
\begin{equation}
Y^2=4X^4 + 2X^3 - 27X^2 - 4X + 61. \label{GPsymquarticec}
\end{equation}

Eq. \eqref{GPsymquarticec} represents the quartic model of an elliptic curve, and the birational transformation defined by
\begin{equation}
\begin{aligned}
X &= -(x - 2y - 7)/(2(4x - 73)), \\
 Y & = (8x^3 - 219x^2 - 4y^2 + 45y + 24187)/(2(4x - 73)^2),
\end{aligned}
\label{GPsymbirat1}
\end{equation} 
and
\begin{equation}
 x = 8X^2 + 2X + 4Y - 9, \quad  y = 32X^3 + 12X^2 + 16XY - 108X + 2Y - 8,
\end{equation}
reduces the quartic curve \eqref{GPsymquarticec} to the Weierstrass form of the elliptic curve given by the cubic equation,
\begin{equation}
y^2 = x^3 - 1227x + 16346. \label{GPsymcubicec}
\end{equation}

A reference to Cremona's  database of elliptic curves \cite{Cr} shows that the curve \eqref{GPsymcubicec} has rank 1 and the single generator of the Mordell-Weil group is the point $R$ whose coordinates are given by  $(x, y) =(25, 36)$. We can accordingly find infinitely many rational points on the curve \eqref{GPsymcubicec} and thus find infinitely many rational solutions of the simultaneous equations \eqref{GPsymcondP2b} and \eqref{GPsymcondP3b}. All of these solutions do not yield the desired examples of geometric progressions on a circle. For instance, the point $R$ and $2R$ do not yield positive values of $r$. However, the point $3R$ with coordinates $(385, 7524)$ yields an example of  the circle,
\begin{equation}
(x - 15888)^2 + y^2 = 19825^2, \label{GPsymcircle}
\end{equation}
on which there are  four points $(3993, 15860), (7623, 18020), (14553, 19780)$, $(27783, 15860)$ whose  abscissae are in geometric progression with the common ratio being $21/11$.

Since the elliptic curve \eqref{GPsymcubicec} has positive rank, it follows from a theorem of Poincar\'e  and Hurwitz \cite[Satz 11, p. 78]{Sk} that  there are infinitely many rational points in the neighbourhood of the rational point $(385, 7524)$ on the  elliptic curve \eqref{GPsymcubicec}, and these rational points would yield infinitely many examples of geometric progressions on a circle with common ratio $> 1$ with the first and fourth  points of the geometric progression being situated symmetrically on either side of the diameter of the circle parallel to the $y$-axis. 

\subsubsection{More general examples}\label{GPgenex}
We will now obtain more general examples in which  the first and last terms of the geometric progression on the circle \eqref{eqcircle} are neither points of intersection of the circle with the $x$-axis, nor are these  points symmetrically situated as in the case of geometric progressions obtained in Section \ref{GPsym}. We will  show how to obtain infinitely many such examples  by finding infinitely many solutions of the simultaneous diophantine equations \eqref{condh12} and \eqref{condh23} using the numerical results given in Table \ref{Tnumex}. 

The solution at S. No. 2 in Table \ref{Tnumex} has $u_2=3/5, u_3= 1/3$, and we substitute these values in the two equations \eqref{condh12} and \eqref{condh23} which may now be written as follows:
\begin{align}
(150r + 48)u_1^2 - 340ru_1 + 150r + 48 & = 0, \label{GPgenequ1}\\
(48u_4^2 + 48)r - 102u_4^2 + 340u_4 - 102 & =0.  \label{GPgenequ4}
\end{align}
We solve Eq. \eqref{GPgenequ1} for $r$ when we get
\begin{equation}
r= -24(u_1^2 + 1)/(5(15u_1^2 - 34u_1 + 15)),\label{valru1}
\end{equation}
and now Eq. \eqref{GPgenequ4} reduces to
\begin{multline}
\quad \quad (4401u_1^2 - 8670u_1 + 4401)u_4^2 - (12750u_1^2 - 28900u_1 + 12750)u_4\\
 + 4401u_1^2 - 8670u_1 + 4401=0. \quad \quad \quad \label{GPgenequ4a}
\end{multline}

Eq.  \eqref{GPgenequ4a} is a quadratic equation in $u_4$ and it will have a rational solution if its discriminant is a perfect square, that is, there must exist rational numbers $u_1$ and $v$ such that the following equation has a rational solution:
\begin{equation}
v^2=1329489u_1^4 - 6745260u_1^3 + 11011078u_1^2 - 6745260u_1 + 1329489. \label{GPgenequ4b}
\end{equation}

Now  Eq. \eqref{GPgenequ4b}  represents the quartic model of an elliptic curve, and the birational transformation defined by
\begin{equation}
\begin{aligned}
u_1 & = (68677x - 2616678847 + 294y)/(252y + 10671x + 1337348919),\\
v & = (3148740x^3 - 701722434420x^2 + 24429873322170120x \\
& \quad \quad - 53230085194000y + 316866956396335870560)\\
& \quad \quad  \times (84y + 3557x + 445782973)^{-2},
\end{aligned}
\end{equation}
and
\begin{equation}
\begin{aligned}
x &= -(10439127u_1^2 - 21636858u_1 - 16065v + 3865183)/(8(6u_1 - 7)^2),\\
y & = -765(17040807u_1^3 + 49967001u_1^2 + 10671u_1v - 171887163u_1\\
& \quad \quad - 68677v + 80372691)/(16(6u_1 - 7)^3),
\end{aligned}
\end{equation}
reduces Eq. \eqref{GPgenequ4b} to the Weierstrass form of the elliptic curve given by
\begin{equation}
y^2 = x^3 + x^2 - 7758767360x + 237867647099508. \label{GPgencubicec}
\end{equation}

The conductor of the elliptic curve \eqref{GPgencubicec} is 43368331440 and the curve is not included in the databases of Cremona and of Stein on elliptic curves. Accordingly, we used the software APECS (a package written in MAPLE) to  determine that the rank of the elliptic curve \eqref{GPgencubicec} is 3, and three independent rational points $R_1, R_2, R_3$ on it are as follows:
\begin{equation*}
\begin{aligned}
R_1 & = (32760151/441, 78279537050/9261),\quad R_2= (64588, 2486862),\\
R_3 & =(5836, 13884750).
\end{aligned}
\end{equation*}

We can now obtain infinitely many rational points on the elliptic curve \eqref{GPgencubicec} using the group law and corresponding to each such rational point, we can obtain a rational solution  of the simultaneous equations \eqref{GPgenequ1} and \eqref{GPgenequ4}. All such solutions will not yield the desired examples of geometric progressions on a circle. For instance, while the point $R_1$ yields the known example of a geometric progression listed at S. No. 2 of Table \ref{Tnumex}, the point $R_2$ leads to a negative value of $r$, and hence we do not get a geometric progression with a positive common ratio. The point $R_3$, however,  yields  a new example of a geometric progression on  the circle 
\begin{equation}
(x-14942502807)^2+y^2=13813414685^2,
\end{equation}
on which the four points,
\begin{align*}
& (632132844180584652/554537833, 303396613510115972/554537833), \\
& (2754195732, 6500430440), (6654453996, 11050731748),\\
 &(1230048832913343556/76505437, 1053225199474661304/76505437),
\end{align*}
have their abscissae  in geometric progression with common ratio $554537833/229516311$.

While all rational points on the elliptic curve \eqref{GPgencubicec} do not yield examples of geometric progressions with a positive common ratio, the aforementioned theorem of Poincar\'e and Hurwitz ensures the existence of infinitely many rational points on the curve \eqref{GPgencubicec} in the neighbourhood of the points $R_1$ and  $R_3$ and these rational points will yield infinitely many examples of geometric progressions on a circle with a positive common ratio and such that the terms of the geometric progression do not contain the points of intersection of the circle with the $x$-axis.

We can similarly use the other solutions listed in Table \ref{Tnumex} to obtain further examples of geometric progressions on a circle  such that the terms of the geometric progression do not contain the points of intersection of the circle with the $x$-axis.

\noindent Ajai Choudhry, 13/4 A Clay Square, Lucknow - 226001, India.

\noindent E-mail address: ajaic203@yahoo.com

\end{document}